\documentclass[10pt,twoside]{article}

\usepackage[centertags]{amsmath}
\usepackage{amsfonts}
\usepackage{amssymb}
\usepackage{amsthm}
\usepackage{epsfig}
\usepackage{color}
\allowdisplaybreaks[4]

\newcommand{\abs}[1]{\lvert #1 \rvert}

\definecolor{c20}{rgb}{0.,0.7,0.}
\definecolor{c30}{rgb}{0.,0.,1.}
\definecolor{c40}{rgb}{1,0.1,0.7}
\definecolor{c50}{rgb}{1,0,0}
\definecolor{c60}{rgb}{1,0.9,0.1}

\def\ff#1{\textcolor{c50}{#1}}
\def\ff#1{#1}

\def\ef#1{\textcolor{c50}{#1}}
\def\ef#1{#1}

\def\aH#1{\textcolor{c30}{#1}}
\def\aH#1{#1}

\def\AE#1{\textcolor{c30}{#1}}
\def\AE#1{#1}
\def\JE#1{\textcolor{c30}{#1}}
\def\jE#1{\textcolor{c30}{#1}}
\def\EEE#1{\textcolor{c30}{#1}}
\def\EEE#1{#1}

\def\vE#1{\textcolor{c30}{#1}}
\def\wE#1{\textcolor{c30}{#1}}
\def\zE#1{\textcolor{c30}{#1}}

\def\xE#1{\textcolor{c30}{#1}}

\def\cE#1{\textcolor{c30}{#1}}

\def\zE#1{#1}
\def\cE#1{#1}

\def\cH#1{\textcolor{c30}{#1}}

\def\yT#1{\textcolor{c50}{#1}}
\def\tT#1{\textcolor{c50}{#1}}
\def\zT#1{\textcolor{c50}{#1}}
\def\zzT#1{\textcolor{c50}{#1}}
\def\cH#1{#1}
\def\xE#1{#1}

\def\wE#1{#1}
\def\yT#1{#1}
\def\vE#1{#1}
\def\tT#1{#1}
\def\jE#1{#1}
\def\JE#1{#1}
\def\zT#1{#1}
\def\zzT#1{#1}

\def\wHH{ \mathcal{H}^2_{2H }}

\newcommand{\COM}[1]{}

\newcommand{\EE}[1]{\mathbb{E}\Bigl\{#1 \Bigr\}}
\newcommand{\Es}[1]{\mathbb{E}\{#1 \}}

\newcommand{\ABs}[1]{ \biggl \lvert #1 \biggr \rvert}
\newcommand{\R}{\mathbb{R}}

\newcommand{\inr}{\in \R}
\newcommand{\BQN}{\begin{eqnarray}}
\newcommand{\EQN}{\end{eqnarray}}
\newcommand{\BQNY}{\begin{eqnarray*}}
\newcommand{\EQNY}{\end{eqnarray*}}
\newtheorem{theo}{Theorem}[section]

\newtheorem{de}[theo]{Definition}
\newtheorem{lem}[theo]{Lemma}

\newcommand{\BL}{\begin{lem}}
\newcommand{\EL}{\end{lem}}

\newcommand{\netheo}[1]{{Theorem \ref{#1}}}
\newcommand{\BT}{\begin{theo}}
\newcommand{\ET}{\end{theo}}
\newcommand{\BD}{\begin{de}}
\newcommand{\ED}{\end{de}}

\def\IF{\infty}
\newcommand{\E}[1]{\mathbb{E}\left \{ #1 \right \}}

\def\ZHT{Z_H(\tau,s)}
\def\ZHTY{Z_H^*(\tau,s)}

\def\MHT{\xE{M_{H}(T)}}

\def\MYA{M_H^*(a,b,T)}

\begin{document}

\centerline{ \Large \bf Large Deviations of Shepp Statistics for Fractional Brownian Motion }

\bigskip
\centerline{Enkelejd Hashorva\footnote{Faculty of Business and Economics,University of Lausanne, Extranef, UNIL-Dorigny, 1015 Lausanne, Switzerland} 
and  Zhongquan Tan\footnote{College of Mathematics, Physics and Information Engineering, Jiaxing University, Jiaxing 314001, PR China}
}

\bigskip
\centerline{\today{}}

{\bf Abstract:} Define  the incremental fractional Brownian field
$\ZHT=B_{H}(s+\tau)-B_{H}(s), H\in (0,1)$, where $B_{H}(s)$ is a
standard fractional Brownian motion with Hurst index
$H\in(0,1)$. In this paper we derive the exact asymptotic
\aH{behaviour} of the maximum
$M_H(T)=\max_{(\tau,s)\in[0,1]\times[0,T]} \ZHT $ \ef{for any  $H\in (0,1/2)$ complimenting thus the result of Zholud
(2008) which establishes the exact tail asymptotic behaviour of
$M_{1/2}(T)$.}

{\bf Key Words:} Shepp statistics; scan statistics; exact asymptotics; extremes; fractional Brownian motion.

{\bf AMS Classification:} Primary 60G15; secondary 60G70

\section{Introduction}
Let $\{B_H(t),t\ge 0\}$ be a standard fractional Brownian motion (fBm) with Hurst index 
$H \in (0,1)$ which is a centered $H$-self-similar Gaussian process with stationary increments, almost surely continuous sample  paths, \aH{$B_H(0)=0$ and its covariance  function is} given by
\BQNY
Cov(B_H(t),B_H(s))=\frac{1}{2}(t^{2H}+s^{2H}-\mid t-s\mid^{2H}),\quad t,s\ge0.
\EQNY
An \JE{important} random field defined in terms of this fBm is  \cE{the so-called incremental fractional Brownian motion}  
$$\ZHT =B_{H}(s+\tau)-B_{H}(s), \quad s,\tau \ge 0. $$
In various statistical applications the incremental fractional Brownian motion  appears as the limit
model. Typically,   when independent and identical observations are modeled, then the limit model has $H=1/2$.
\JE{A closely related random field, namely the standardised \aH{incremental fractional Brownian motion}
$$\ZHTY =\frac{B_{H}(s+\tau)-B_{H}(s)}{\tau^H}, \quad s,\tau \in (0,\IF) $$
serves \ff{also} in various applications as a limit model.} For instance with motivation from queuing
theory, consider $\{K(t), t\ge 0\}$ a homogeneous Poisson process
with $\Es{K(t)}=\lambda t, \lambda>0$ and set for some $T,\tau$
positive $K(\tau,T):= \sup_{0 \le s \le T} (K(s+\tau)-K(s))$. The
random variable $K(\tau,T)$ is the maximum service length of an
$M/G/\IF$ queue with deterministic service time $\tau$; it is also the version of
scan statistics on the positive half line (see e.g., Cressie
(1980)).  For the study of  $K(\tau,T)$ the following  convergence
in distribution
 $$ \frac{ K(\tau,T)- \lambda \tau}{\sqrt{\lambda \ff{\tau}}}
 \to \sup_{0 \le s \le T} \ef{Z_{1/2}}^*(\tau,s)
 , \quad \lambda \to \IF$$
is important since the distribution function of $\sup_{0 \le s \le T} \ef{Z_{1/2}}^*(\tau,s)$ is derived in Shepp (1971),
see also Slepian (1961), Shepp (1966) and Theorem 3.2 in Cressie (1980).\\
Various authors refer to the process $\{Z_H^*(\tau,T), \tau\ge 0\}$
as the standardised Shepp statistics. Important results for Shepp
statistics  and related quantities can be found in Deheuvels and
Devroye (1987),
Siegmund and Venkatraman (1995), D\"umbgen and Spokoiny (2001), Kabluchko and Munk (2008) and Zholud (2009). \\
The recent papers Zholud (2008) and Kabluchko (2007,2011a) present asymptotic results on the extremes of Shepp statistics and \JE{standardised
Shepp statistics},  i.e., therein the tail asymptotic behaviour of
$$ \MHT = \sup_{\tau \in [0,1], 0 \le s \le T} Z_H(\tau,s) \quad \text{ and }  \MYA= \sup_{\tau \in [a,b], 0 \le s \le T}Z_H^*(\tau,s),
\quad 0 < a < b < \IF$$
 for the case $H=1/2$ \JE{is investigated} dealing thus with the increments of the Brownian motion.

\jE{In view of Zholud (2008)
for any  $T>0$
\begin{eqnarray}\label{ZH08}
P\left(\zT{M_{1/2}(T)} >u\right) =\widetilde{\mathcal{H}}_* T u^{2}\Psi(u)(1+o(1)), \quad
u\rightarrow\infty,
\end{eqnarray}
where $\Psi$ is the survival function of a $N(0,1)$ random variable and \aH{the constant $\widetilde{\mathcal{H}}_*$ is given by}
$$\widetilde{\mathcal{H}}_*=\lim_{\cE{a}\rightarrow\infty}\lim_{\cE{b}\rightarrow\infty}a^{-1}e^{-\frac{a+b}{2}} \cE{\mathbb{E}}
\biggl\{\exp\bigg(\max_{0\leq t\leq a \atop 0\leq s\leq
b}B_{1/2}(t+s+a)-B_{1/2}(t)\bigg)\biggr\}.$$ For \aH{any Hurst
index} $H\not=1/2$ the independence of the increments of $B_H$ does
not hold, which has been the crucial  property in the derivation of
\eqref{ZH08}. In our main result given in
\netheo{Th:main} we derive the exact asymptotics of $M_H(T)$ for $H \in (0,1/2)$, \ef{which is the well-known short-range dependence case for fBm}. If $H\in (1/2,1)$, \ef{thus we have a long-range dependence,}  we have a \aH{much} more involved problem which will be therefore considered elsewhere. \\ \ff{Numerous authors have considered properties and characterisations of fBm, see e.g., Mishura and Valkeila (2011), Kabluchko (2011b) and the references therein.
Our contribution present a new result for the Shepp statistics of fBm, which we believe is important for both future theoretical and applied developments.}\\
Clearly, our result  on the tail asymptotic behaviour of $M_H(T)$ implies certain asymptotic bounds for the tail asymptotics
of $\MYA$.}  The exact asymptotics of $\MYA$ is however easier to deal with and follows by a direct application of the results of Chan and Lai  (2006) since the pertaining random field is locally stationary; see Mikhaleva and Piterbarg (1996), and Piterbarg (1996) for the main findings concerning locally stationary random fields.

Brief outline of the rest of the paper: Section 2 displays the main result, its proof is given in Section 3.

\section{Main Result}
 In the asymptotic theory of Gaussian processes two important constants are crucial, namely the Pickands and Piterbarg constants.
 Since in our results only the former constant appears, we briefly mention that it is defined by  (see Pickands (1969), Piterbarg (1996))
  $$\mathcal{H}_{2H}=\lim_{\lambda\rightarrow\infty} \lambda^{-1} \E{ \exp\left(\max_{t\in[0,\lambda]}\ef{\Bigl(}\sqrt{2}B_{H}(t)-t^{2H}\ef{\Bigr)}\right)} \in (0,\IF).$$
  It is well-known that $\mathcal{H}_{1}=1$ and $\mathcal{H}_{2}= 1/ \sqrt{ \pi}$;
   \cE{see Piterbarg (1972) which gives the first rigorous proof of Pickands theorem presented in Pickands (1969),
   D\c{e}bicki (2002), Wu (2007), D\c{e}bicki and Kisowski  \aH{(2009)} and \JE{D\c{e}bicki and Tabi\'{s} (2011)}  for generalisations of Pickands constant.}

The result of (\ref{ZH08}) is of some importance for dealing with the general case $H\not=1/2$. However
 we cannot use the method of proof \ff{in} Zholud (2008) which relies on the independence of increments of Brownian motion.
\aH{Our proof of the main result presented below is strongly motivated by the method utilised in the seminal contribution  Piterbarg (2001).}

\BT\label{Th:main}  \aH{For any $H\in(0, 1/2)$} and any $T>0$
\begin{eqnarray}\label{f1}
&&P\left(\MHT >u\right)
=\frac{T}{H} \left(\frac{1}{2}\right)^{1/H}\wHH  u^{\frac{2}{H}-2}\Psi(u)(1+o(1))
\end{eqnarray}
\aH{holds as $u\to \IF.$}
\ET

\textbf{Remark}.  a) \ef{As in Zholud (2008) also for the case $H \ff{\in (0,1/2)}$ it is possible to derive a similar expansion as in \eqref{f1} when $T=T_u$  depends} on $u$. From our proof we see that
the result holds if $T=T_u$ is such \EEE{\ff{that} $\lim_{u \to \IF} T_u u^{1/d}=\IF$} for some \EEE{$d\in(\zE{H,1/2})$}. \aH{Following word-by-word the} arguments of H\"usler and Piterbarg (2004) the proof of  the case  $ T_u<\exp(cu^{2})$ can also be included when  $c\in(0,1/2)$.

b) \ef{Borrowing the arguments of} \aH{H\"{u}sler and Piterbarg (2004),  we \ff{obtain that}} for $H\in (0,1/2)$ the Gumbel limit law
\BQN\label{Gumbel} \lim_{T \to \IF}  \xE{\max}_{x\inr } \ABs{
P\left(a_{T}( \MHT  -b_{T})\leq x\right)-\exp(-e^{-x})}=0, \EQN
\ff{holds}, where
\BQN \label{abT}
a_{T}=\sqrt{2 \ln T},\ \ b_{T}=a_T+a_T^{-1} \Bigl[ (\frac{1}{H}-\frac{3}{2}) \ln  \ln T+ \ln (2^{-3/2} \wHH H^{-1}(2\pi)^{-1/2})\Bigr].
\EQN


\section{Proofs}

\aH{We \ff{present} first a lemma which is crucial for the proof of Theorem 2.1.}

\def\ve{\varepsilon}
\def\DU{\xE{\delta_u}}

\def\SS{\ln^2 u}
\BL \label{Lemma 4.1}
Let $\ve \in (0,H), H\in (0,1)$ and let \ff{$T>0$} be given. Then for all large $u$ and $\DU:= \ln^2 u /u^2$
\begin{eqnarray}
\label{eq4.1}
P\left( \xE{\max}_{\tau\in[0, 1- \DU]\atop s\in [0,T]} \ZHT \xE{> u} \right)\leq CTu^{2/H-1}\exp\left(-\frac{1}{2}u^{2}-\aH{(H-\ve)} \ln ^{2}u\right)
\end{eqnarray}
holds with $C>0$ not depending on $T$ and $u$.
\EL

\textbf{Proof:} \xE{For any $\tau,s, \tau',s'$ positive we have
$$ \EE{ (\ZHT- Z_H(\tau',s'))^2} \le G [\abs{\tau -\tau'}^{2H}+\abs{s-s'}^{2H}],$$
\zzT{with some constant $G>0$}. Consequently, by Piterbarg
inequality (given in \ff{Theorem 8.1 in Piterbarg (1996), see also Theorem 8.1 in Piterbarg (2001)} and Proposition 3.2 in Tan and Hashorva (2013)) for any $x \in (0,1)$ we have
\BQNY P\left( \xE{\max}_{\tau\in[0,x]\atop s\in [0,T]} \ZHT > u
\right) \leq
CT u^{2/H-1}\exp\left(-\frac{u^{2}}{2x^{2 H}}\right)
\EQNY
for some $C$ independent of $x$ and $u$. For $x=1- \SS/u^{2}$ we obtain
\BQNY
P\left( \xE{\max}_{\tau\in[0,x]\atop s\in [0,T]} \ZHT > u \right)
&\leq& CTu^{2/H-1}\exp\left(-\frac{u^{2}}{2(1- \SS/u^2)^{2H}}\right),
\EQNY
hence the proof follows.} \hfill$\Box$


\def\TE{\AE{T_\varepsilon}}
\def\Su{S_u}
\def\upH{\Su^{-1}}
\def\umH{\jE{u^{\frac{1}{2H'}}}}
\def\umH{\Su}
\def\DU{\xE{\delta_u}}

\def\SS{\ln^2 u}

\def\TE{\AE{T_\varepsilon}}
\def\Su{S_u}
\def\upH{\Su^{-1}}
\def\upH{R_u}
\def\umH{\jE{u^{\frac{1}{2H'}}}}
\def\umH{{R_u^{-1}}}

\textbf{Proof of Theorem 2.1.} \EEE{Let in the following $\upH=u^{-2/H'}$ for some $H'\in (2H,1)$ and let $\DU= u^{-2}\ln^2 u $.
By Lemma 3.1, we can restrict the considered domain
of $(\tau,s)$ to  $\tau\in[1- \DU, 1]$.
For this choice of $\upH$, 
we have
\BQN \label{DSu}
\lim_{u\to \IF}   \frac{\DU}{\upH}= \IF, \quad  
\text{  and  }\lim_{u\to \IF}   u^2 \upH= 0.
\EQN
Note that with $\mathcal{A}:= \{(t,s) \inr^2: t\in[0,T+1],s\in [0,T], t-s\in[0,1]\}$
$$M_{H}(T)=\max_{(t,s) \in \mathcal{A}}Y_{H}(t,s), \quad \text{ where  }Y_{H}(t,s)=B_{H}(t)-B_{H}(s).$$}
Define next for any $l=1,\cdots, [T \umH ]$
$$J_{k,l}=[1+(l -k) \upH ,
1+(l -k+1) \upH ]\times[(l-1) \upH ,
l \upH ],$$
with $k=1,\cdots, [\DU  \umH +1]$ and
$$J'_{k,l}=[1+(l-1- k)\upH ,
1+(l-k) \upH ]\times[(l-1) \upH ,
l \upH ],$$
where $k=1,\cdots, [\DU  \umH -1]$.\\
The variance function $\vE{\sigma^{2}(\tau,s)}$ of $ \ZHT $ equals
$\tau^{2H}$. \ef{Consequently, $\vE{\sigma}(1,s)=1$ for all $s \in [0,1]$ and hence the maximum point of the variance is not a single point but taken on $(1,s), s\in (0,1)$.} 
Taylor \wE{expansion yields}  \BQNY
\sigma(\tau,s)=\sigma_{Z}(\tau)=1-H(1-\tau)+o(1-\tau)^{2} \EQNY as
$\tau\uparrow1$. \COM{ By series expansion we find for any
$\tau,\tau'$ with $0<\tau_{1}<\tau,\tau'<\tau_{2}<1$
\BQNY
r(\tau,s;\tau',s')\leq C|s-s'|^{2H-2}
\EQNY
for some constant $C>0$ and all $s,s'$ with $|s-s'|$ sufficiently large}
 \yT{Hence, for arbitrarily small $\varepsilon>0$ and for all sufficiently large $u$,  on $ \jE{J_{k,l}} $ the variance
$\sigma^{2}_{Y}(t,s)$ of $Y_{H}(t,s) $ satisfies
 $$\COM{1-(H+\varepsilon)k \upH \leq }\sigma_{Y}(t,s)\leq 1-(H-\varepsilon)(k-2) \upH .$$
On $ \jE{J'_{k,l}} $ the variance $\sigma_{Y}^{2}(t,s)$ of $Y_{H}(t,s) $
satisfies
 $$1-(H+\varepsilon)(k+1) \upH \leq \sigma_{Y}(t,s)\COM{\leq
 1-(H-\varepsilon)(k-1) \upH }.$$}
 \yT{For the  correlation function
$r_{Y}(t,s;t',s')$ of $Y_{H}(t,s) $, we have
\BQNY
r_{Y}(t,s;t',s')=1-\frac{1}{2}(1+o(1))(|t-t'|^{2H}+|s-s'|^{2H})
\EQNY
as $t-s,t'-s'\uparrow1$, $s-s'\rightarrow 0$ and $t-t'\rightarrow
0$.} Consider the \ef{centered and homogeneous}
Gaussian \cH{random} field $\zeta_{\wE{H}\pm\varepsilon}(t,s)$ with
covariance function
$$\exp\left(-(\frac{1}{2}\pm\varepsilon)(|s-s'|^{2H}+|t-t'|^{2H})\right), \quad s,s',t,t'\inr,$$
\aH{where $\varepsilon\in (0,H)$.}  \ef{From this point, we assume below that $H\in (0, 1/2)$.} Using Slepian inequality and
Theorem 7.2. of Piterbarg (1996) for all sufficiently large $u$ we
have
\begin{eqnarray}
\label{eqUP}
\yT{P\left( \xE{\max}_{(t,s)\in J_{k,l}} Y_{H}(t,s)  >u\right)\leq   {\wE{(\frac{1}{2}+\varepsilon)^{\frac{1}{H}}}\wHH u_{k-}^{\frac{2}{H}} \upH^{2} \Psi(u_{k-})}=:\eta(\varepsilon,u_{k-})}
\end{eqnarray}
and
\begin{eqnarray}
\label{eqLO}\yT{P\left( \xE{\max}_{(t,s)\in J'_{k,l}} Y_{H}(t,s)  >u\right)\geq
{\wE{(\frac{1}{2}-\varepsilon)^{\frac{1}{H}}}\wHH u_{k+}^{\frac{2}{H}} \upH^{2}\Psi(u_{k+})}=:\eta(-\varepsilon,u_{k+})},
\end{eqnarray}
where
$$u_{k-}=\frac{u}{1-(H-\varepsilon)(k-2) \upH },\ \ u_{k+}=\frac{u}{1-(H+\varepsilon)(k+1) \upH }.$$
\def\HE{\xE{H_\varepsilon} }
Further, for all sufficiently large $u$  we \ff{obtain utilising} further \eqref{DSu}  
(set $\HE:=
H - \varepsilon, \TE:= (1/2+ \varepsilon)^{\frac{1}{H}} T \wHH$)
\begin{eqnarray}
\label{eqp2.1}
P\left( \xE{\max}_{(t,s) \in \mathcal{A} }  Y_{H}(t,s)  >u\right)
&\leq& P\left( \xE{\max}_{(t,s)\in \cup J_{k,l}}  Y_{H}(t,s)  >u\right)\nonumber\\
&\leq& \sum_{k=1}^{[\umH \DU   +1]}\sum_{l=1}^{[T \umH ]+1}\cE{\eta(\varepsilon,u_{k-})}\notag \\
&=&([T \umH ]+1)\sum_{k=1}^{[\umH \DU   +1]}
\cE{\eta(\varepsilon,u_{k-})}\notag\\
&=&\aH{\TE \frac{1+o(1)}{\sqrt{2\pi}}u^{\frac{2}{H}-1} \jE{e^{- \frac{u^2}{2}}}
\sum_{k=1}^{[ \umH   \DU +1]}\frac{e^{-\frac{1}{2}(u_{k-}^{2}-\jE{u^2})} \upH }
{(1-\HE (k-2) \upH )^{\frac{2}{H}-1}}}\nonumber\\
&\le&\aH{\frac{\TE}{1-\varepsilon'} \xE{u^{\frac{2}{H}} \Psi(u)}
\sum_{k=1}^{[ \umH   \DU ]+2} \exp\Bigl( \frac{-u^2}{2}\Bigl(  \frac{1}{ (1- H_\ve(k-2) \upH)^2}-1\Bigr) \Bigr) \upH}
\nonumber\\
&\jE{\le}&\frac{\TE}{1-\varepsilon'} \xE{u^{\frac{2}{H}- \jE{2}} \Psi(u)}
\sum_{k=1}^{[ \umH  \DU]+2} \exp\Bigl( -\frac{\HE}{1+\varepsilon'''}(k-2)u^{2} \upH   \Bigr) u^{2} \upH
\nonumber\\
&\jE{=}&\frac{\TE}{1-\varepsilon'} \xE{u^{\frac{2}{H}- \jE{2}} \Psi(u)}
\int_0^{\aH{\IF} }  \exp\Bigl(- \frac{\HE}{1+\varepsilon'''} x\Bigl)\, dx (1+o(1))\notag\\
&=&\frac{\TE (1+\varepsilon''')}{ \AE{\HE}(1- \varepsilon')}u^{\frac{2}{H}-2}\AE{\Psi(u)}(1+o(1)), \quad u\to \IF, 
\end{eqnarray}
where $\varepsilon',\varepsilon'', \varepsilon''' \in (0,1)$ above are appropriately chosen constants and the passing from the sum to
the integral is legitime since \eqref{DSu} holds.
\COM{
\begin{eqnarray*}
&&T(\frac{1}{2}+\varepsilon)\wHH \xE{u^{\frac{2}{H}-2} \Psi(u)}
\int_{0}^{u\ln u}\frac{1}{(1-\xE{\HE x/u^{2}})^{\frac{2}{H}-\frac{2}{H'}-1}}\exp\left(-
\frac{1}{2} \frac{ 2 \HE x - (\HE x/u^{2})^2}{(1-\HE x/u^2)^{2}}\right)
dx(1+o(1))\nonumber\\
&&=T(\frac{1}{2}+\varepsilon)\wHH \xE{u^{\frac{2}{H}-2} \Psi(u)}
\xE{\int_{0}^{\infty}\exp(- \HE x) dx} (1+o(1))\nonumber\\
&&=T(\frac{1}{2}+\varepsilon)\wHH u^{\frac{2}{H}-2}
\frac{1}{\HE }\Psi(u)(1+o(1)).
\end{eqnarray*}
}
Similarly, for all sufficiently large $u$
\begin{eqnarray}
\label{eqp2.2}
&&P\left( \xE{\max}_{(t,s)\in \mathcal{A}}  Y_{H}(t,s)  >u\right)\nonumber\\
&&\geq P\left( \xE{\max}_{(t,s)\in \cup J'_{k,l}}  Y_{H}(t,s)  >u\right)\nonumber\\
&&\geq \sum_{l=1}^{[T \umH ]}\sum_{k=1}^{[ \umH   \DU -1]}
\cE{\eta(-\varepsilon,u_{k+})}
-\sum_{(l,k), (l',k') :\rho((s,t),(s',t'))\geq  \upH}P\left( \xE{\max}_{(t,s)\in  \jE{J'_{k,l}} }  Y_{H}(t,s) >u,
 \xE{\max}_{(t',s')\in  \jE{J'_{k',l'}} } Y_{H}(t',s')  >u\right)\nonumber\\
&&\ \ \ -\sum_{(l,k),(l',k'):\rho((s,t),(s',t'))< \upH \atop (l,k)\neq(l',k')}P\left( \xE{\max}_{(t,s)\in  \jE{J'_{k,l}} }  Y_{H}(t,s)  >u,
 \xE{\max}_{(t',s')\in  \jE{J'_{k',l'}} }  Y_{H}(t',s')  >u\right),
\end{eqnarray}
where \ff{$\rho(\cdot, \cdot)$} is \aH{the Euclidean distance of two points in $\R^2$}. The \cE{first} sum \cE{can be bounded from below with}
the same arguments as in the proof of (\ref{eqp2.1}) by
$$\vE{(1/2-\varepsilon)^{\frac{1}{H}}}\frac{T}{H+\varepsilon} \wHH u^{\frac{2}{H}-2}
\Psi(u)(1+o(1)).$$ \tT{A generic term of the second sum can be estimated as follows. The Gaussian field
$ W_H(t,s,t',s')= Y_{H}(t,s)+Y_{H}(t',s')$ on $ \jE{J'_{k,l}} \times  \jE{J'_{k',l'}} $ \aH{has for $H \in (0,1/2)$ variance function}
\begin{eqnarray*}
2+2r_{Y}(t,s;t',s')&=&4-(|t-t'|^{2H}+|s-s'|^{2H})(1+o(1)) \leq 4-\EEE{(2- \varepsilon)}R_u^{2H},
\end{eqnarray*}
\aH{where the inequality holds for all large $u$ and $\ve\in (0,H)$.} Consequently, for some $C>0,c_*>0$ by Piterbarg inequality
\begin{eqnarray*}
P\left( \xE{\max}_{(t,s)\in  \jE{J'_{k,l}} }  Y_{H}(t,s) >u,
 \xE{\max}_{(t',s')\in  \jE{J'_{k',l'}} } Y_{H}(t',s')  >u\right)
&\leq & P\left( \xE{\max}_{(t,s, t',s')\in  \jE{J'_{k,l}} \times  \jE{J'_{k',l'}} }  W_H(t,s,t',s') >2u\right)\\
&\leq & C \EEE{R_u^2}u^{\frac{4}{H}}\Psi\left(\frac{u}{\sqrt{1-\EEE{(1/2- \varepsilon)R_u^{2H}}}}\right)\\
&=& C \EEE{R_u^2}u^{\frac{4}{H}} \exp\left(-c_*u^{2} R_u^{2H}\right)\Psi(u)\\
&=& o( u^a \Psi(u)), \quad u\to \IF
\end{eqnarray*}
for any real $a$, \aH{since $H< 1/2$ and $H' \in (2H,1)$  imply
$$u^2 R_u^{2H}= u^{2(1 - 2H/H')} \to \IF, \quad u\to \IF.$$
}
Therefore, in order to complete the proof we need to consider only} the third term on the right-hand side of (\ref{eqp2.2}). The generic term in this sum is equal
to{
\begin{eqnarray}
\label{eqp2.3}
&&P\left( \xE{\max}_{(t,s)\in  \jE{J'_{k,l}} }  Y_{H}(t,s) >u\right)+P\left( \xE{\max}_{(t,s)\in  \jE{J'_{k',l'}} }  Y_{H}(t,s) >u\right)
-P\left( \xE{\max}_{(t,s)\in  \jE{J'_{k,l}} \cup  \jE{J'_{k',l'}} }  Y_{H}(t,s) >u\right)\nonumber\\
&&=2P\left( \xE{\max}_{(t,s)\in  \jE{J'_{k,l}} }  Y_{H}(t,s) >u\right)-P\left( \xE{\max}_{(t,s)\in  \jE{J'_{k,l}} \cup  \jE{J'_{k',l'}} }  Y_{H}(t,s) >u\right).
\end{eqnarray}
The first term on the right-hand side of (\ref{eqp2.3}) has been
already estimated. For the second one  note that
$\rho(s,t;s',t')< \upH $, so both $l-l'$ and $k-k'$
cannot be greater than 2. The variance of the \ff{random} field $\{Y_{H}(t,s)$,
$(t,s)\in  \jE{J'_{k,l}} \cup  \jE{J'_{k',l'}} \ff{\}}$, satisfies
 $$1-(H+\varepsilon)(k+1) \upH \leq \sigma_{Y}(t,s)\leq
 1-(H-\varepsilon)(k-2) \upH .$$
Hence we use (\ref{eqUP}) and (\ref{eqLO}) with the same definition
of $u_{k\pm}$. \ff{Consequently}, summing (\ref{eqp2.3}) we obtain a sum which
tends to an integral and which terms are uniformly negligible with
respect to the corresponding terms in the first sum in the
right-hand side of (\ref{eqp2.2}).  \aH{Finally, letting $\varepsilon+ \ve'+\ve''+\ve''' \downarrow0$ establishes the proof. }  \hfill$\Box$

{\bf Acknowledgement}: We are thankful to the referees of the paper as well as to Krzysztof D\c{e}bicki, Lanpeng Ji, Yuliya Mishura and Yimin Xiao 
for numerous suggestions which improved this paper significantly. 
E. Hashorva kindly acknowledge support by the Swiss National Science Foundation Grants 200021-1401633/1 and 200021-134785 as well as 
by the project RARE -318984, a Marie Curie International Research Staff Exchange Scheme Fellowship within the 7th European Community Framework Programme. Tan's work was supported by National Science Foundation of China (No. 11071182) and Research Start-up Foundation of Jiaxing University 
 (No. 70512021).

\end{document}